\documentclass[12pt]{article} 
\usepackage{amsmath, amssymb} 
\usepackage{amsthm} 
\usepackage{enumerate, url} 
\usepackage{graphicx}

\newtheoremstyle{plainsl}%
	{\topsep}
	{\topsep}
	{\slshape} 
	{}
	{\normalfont\bfseries}
	{.}
	{ }
	{}

\swapnumbers

{\theoremstyle{plainsl}
\newtheorem{theorem}{Theorem}[section]
\newtheorem{lemma}[theorem]{Lemma}
\newtheorem{corollary}[theorem]{Corollary}}
{\theoremstyle{remark}
}

\newcommand\lref[1]{Lemma~\ref{lem:#1}}
\newcommand\tref[1]{Theorem~\ref{thm:#1}}
\newcommand\cref[1]{Corollary~\ref{cor:#1}}
\newcommand\sref[1]{Section~\ref{sec:#1}}

\renewcommand\proof{\noindent\textsl{Proof. }}
\newcommand\sqr[2]{{\vbox{\hrule height.#2pt
    \hbox{\vrule width.#2pt height#1pt \kern#1pt
        \vrule width.#2pt}\hrule height.#2pt}}}
\renewcommand\qed{%
	\ifmmode\eqno\sqr53
	\else\nolinebreak\ \hfill\sqr53\medbreak\fi}

\numberwithin{equation}{section}

\newcommand\al{\alpha}
\newcommand\be{\beta}

\newcommand\ga{\gamma}
\newcommand\Ga{\Gamma}

\newcommand\cA{{\mathcal A}}
\newcommand\cB{{\mathcal B}}

\newcommand\cx{{\mathbb C}}

\newcommand\ints{{\mathbb Z}}
\newcommand\re{{\mathbb R}}

\newcommand\sbs{\subseteq}

\newcommand\seq[3]{#1_{#2},\ldots,#1_{#3}}

\DeclareMathOperator{\sym}{Sym}
\newcommand\pmat[1]{\begin{pmatrix} #1 \end{pmatrix}}
\newcommand\sm[3]{\sum_{#1=#2}^{#3}}

\newcommand\ass{\cA}
\newcommand\assb{\cB}
\newcommand\clc{{\cal C}_4}
\newcommand\cxa{\cx[\ass]}
\newcommand\dcfr{{\cal D}_4}
\newcommand\pn{^{\otimes n}}
\newcommand\hn[1]{H(n,#1)}
\newcommand\hnass{\hn{\ass}}
\newcommand\hnq{\hn{q}}

\title{Generalized Hamming Schemes} 

\author{
	Chris Godsil\\
	Combinatorics \& Optimization\\
	University of Waterloo}

\begin{document}
\maketitle
	
\begin{abstract}
	We introduce a class of association schemes that generalizes the
    Hamming scheme.  We derive generating functions for their eigenvalues,
    and use these to obtain a version of MacWilliams theorem.
\end{abstract}

\section{Association Schemes}

An association scheme is a collection of graphs that fit together nicely.
It is easiest to offer a precise definition in terms of adjacency matrices.
So an \textsl{association scheme} $\cA$ on $v$ vertices is a set of
$v\times v$ $01$-matrices $\seq A1d$ such that:
\begin{enumerate}[(1)]
    \item 
    $A_0=I$.
    \item
    $\sm i0d A_i=J$.
    \item
    If $A_i\in\cA$ then $A_i^T\in\cA$.
    \item
    There are constants $p_{i,j}(k)$ such that
    \[
        A_iA_j =\sm r0d p_{i,j}(r)A_r.
    \]
    \item
    For all $i$ and $j$ we have $A_iA_j=A_jA_i$.
\end{enumerate}
By virtue of (4), the complex span of the matrices $A_i$
forms an algebra, known as the \textsl{Bose-Mesner algebra} of 
the scheme and denoted by $\cxa$. Given (5), this algebra is commutative. 
Note that the constants $p_{i,j}(k)$ must be integers; they are the
\textsl{intersection numbers} of the scheme.

An association scheme is \textsl{symmetric} if each matrix $A_i$ is symmetric.
If we restrict ourselves to symmetric schemes the last axiom is redundant.
Some workers prefer to drop the last axiom but since there are as yet no
really convincing combinatorial applications making use of this generality,
we will work with the axioms just stated. The schemes we work with in this
paper will be symmetric.

We view the matrices $\seq A1d$ as the adjacency matrices of directed
graphs $\seq X1d$ and we will often view a scheme as a set of graphs.
We refer to these graphs as the \textsl{color classes} of the scheme.
Since $J$ lies in the Bose-Mesner algebra (by (2)) and since the Bose-Mesner
algebra is commutative, we have $A_iJ=JA_i$ for each $i$. Hence each of the
color classes is regular. 

For more background on association schemes see \cite{delphd,bcn}.

The \textsl{Hamming scheme} $H(n,q)$ is a useful and relevant example of an 
association scheme. Its vertices are the $n$-tuples of elements of a set
$Q$ with size $q$. Two tuples $\al$ and $\be$ are adjacent in the graph $X_i$
if they differ in exactly $i$ positions---equivalently if they are at Hamming
distance $i$. In this context we often call the vertices \textsl{words}.
It is easy to see that the matrices $A_i$ satisfy the first three axioms
and that they are symmetric. We leave the verification of the final two
axioms as an exercise and note that $H(n,q)$ is a symmetric scheme.
The Hamming scheme is of fundamental importance in coding theory.

We introduce a generalization of the Hamming scheme.  Let $\ass$ be a
fixed association scheme with $d$ classes and vertex set $V$.  If $v$
and $w$ are elements of $V^n$, let $h(v,w)$ be the vector of length
$d+1$ with $r$th-entry equal to the number of coordinates $j$ such
that $v_j$ and $w_j$ are $r$-related in $\ass$.  For any $n$-tuples
$v$ and $w$ the vector $h(v,w)$ has non-negative integer entries, and
these entries sum to $n$.  Conversely, any such vector can be written
as $h(v,w)$ for some $v$ and $w$.  If $x$ is an integer vector of
length $d+1$ with entries summing to $n$, let $A_x$ be the $01$-matrix
with rows and columns indexed by $V^n$, and with $A_{v,w}$ equal to 1
if and only if $h(v,w)=x$.  Denote this set of matrices by $\hnass$.

Two things should become clear without too much effort: if
$\ass=H(1,q)$ then $\hnass$ is the Hamming scheme $H(n,q)$ and, in any
case, these matrices in $\hnass$ satisfy the first three axioms for an
association scheme. Generalized Hamming schemes were defined by Delsarte
\cite[Section~2.5]{delphd}, where they are referred to as \textsl{extensions}
and a number of examples are given. One of the goals of this paper is prove that,
for all choices of base scheme, $\hnass$ is an association scheme
and to provide some of its applications. Delsarte's construction can be viewed
as based on the symmetric group $\sym(n)$, we will see that it can be generalized
by replacing the symmetric group by any permutation group on $n$ symbols.
Some of our results have already been applied by Martin and Stinson
in their work on $(T,M,S)$-nets \cite{MarSt} and by Chan and Munemasa
in \cite{ACAM}. (They cite the research report CORR98-20, issued by the 
Department of Combinatorics and Optimization at the University of Waterloo.)

\section{Idempotents and Eigenvalues}

If $M$ and $N$ are $m\times n$ matrices, their \textsl{Schur product} is
the $m\times n$ matrix $M\circ N$ given by
\[
    (M\circ N)_{i,j} := M_{i,j}N_{i,j}.
\]
A matrix is idempotent with respect to the Schur product if and only
if it is a $01$-matrix. The $m\times n$ matrices over a field form 
a commutative algebra relative to the Schur product with multiplicative
identity the all-ones matrix $J$. It is an easy exercise to see that
that a vector space of matrices is an algebra relative to the Schur
product if and only if it has a basis of $01$-matrices. Consequently
the Bose-Mesner algebra of an association scheme is Schur-closed.
It is a crucially important property of an association scheme that
the operations of matrix and Schur multiplication are dual to each other.

An idempotent in algebra is \textsl{minimal} if there is non non-trivial
way of writing it as a sum of other idempotents. The matrices $A_i$
are the minimal Schur idempotents in the Bose-Mesner algebra, viewed as
an algebra under Schur multiplication.
They form one basis of the Bose-Mesner algebra.
We introduce a second important basis. An eigenspace of $\ass$
is a subspace $U$ of $\re^n$, maximal subject to the condition that
each matrix in $\ass$ acts on $U$ as a scalar.  As the Bose-Mesner
algebra is closed under Hermitian transposes, it is semisimple.  It
follows that there are orthogonal projections $\seq E0d$ and real
numbers $p_i(j)$ such that:
\begin{enumerate}[(a)]
    \item
    $E_0={1\over v}J$,
    \item
    $\sum_i E_i=I$,
    \item
    $E_i^2=E_i$ and $E_iE_j=0$ if $i\ne j$,
    \item
    $A_i=\sum_j p_i(j)E_j$. 
\end{enumerate}
We will call these matrices the \textsl{principal idempotents} of the scheme.
(They are the minimal idempotents of the Bose-Mesner algebra, relative
to the usual matrix multiplication.)

The scalars $p_i(j)$ are the \textsl{eigenvalues} of the association
scheme and the matrix $P$ given by
\[
    P_{i,j} =p_j(i)
\]
is its \textsl{matrix of eigenvalues}.  We have
\[
    A_i= \sum_j p_i(j)E_j.
\]
Since the principal idempotents
of $\ass$ lie in its Bose-Mesner algebra, there are also scalars
$q_i(j)$ such that
\[
    E_i ={1\over v}\sum q_i(j)A_j.
\]
If $Q$ is the matrix with $ij$-entry $q_j(i)$ then $PQ=vI$.

Since the Bose-Mesner algebra of a scheme is closed under Schur multiplication,
it follows that there must be scalars $q_{i,j}(k)$ such that 
\[
    E_i\circ E_j = \frac1v \sm r0d q_{i,j}(r)E_r.
\]
These scalars are the \textsl{Krein parameters} of the scheme, they need
not be integers, but must be non-negative.

\section{Automorphisms and Subschemes}

We call an association scheme $\assb$ a \textsl{subscheme} of a scheme
$\ass$ if each matrix in $\assb$ is a sum of elements of $\ass$.
(This is also called a fusion scheme.)  It is immediate that
$\assb$ is a subscheme of $\ass$ if and only the Bose-Mesner algebra
of $\assb$ is a subalgebra of the Bose-Mesner algebra of $\ass$ that
is closed under Hermitian transpose and the Schur product.
We also point out here that a commutative algebra of real matrices
is the Bose-Mesner algebra of an association scheme if and only if
it is closed under transpose and Schur multiplication and contains $J$.

Let $\ass$ be an association scheme with $d$ classes.  Let $M^*$
denote the conjugate-transpose of the complex matrix $M$.  We call a
map $M\mapsto M^\psi$ on the Bose-Mesner algebra of $\ass$ an
\textsl{algebra automorphism} if
\begin{align*}
    (M^*)^\psi &=(M^\psi)^*,\\
    (MN)^\psi &= M^\psi N^\psi,\\
    (M\circ N)^\psi &= M^\psi\circ N^\psi.
\end{align*}
An algebra automorphism must permute the matrices $A_i$ among
themselves and fix $I$.  Similarly it must permute the principal
idempotents $E_j$ and fix $J$.

If $\ass$ be an association scheme that is not symmetric, then the
transpose map is an algebra automorphism.  More important, for us, is
the next example.  Consider the product scheme $\ass\pn$ and let $\Ga$
be any permutation group on the set $\{1,2,\ldots,n\}$.  If $g\in\Ga$
define
$$
(A_{i_1}\otimes\cdots\otimes A_{i_n})^g 
	=A_{i_{1g}}\otimes\cdots\otimes A_{i_{ng}}.
$$
It is immediate that $g$ is an algebra automorphism of $\ass\pn$.

\begin{theorem}\label{thm:indsub}
Let $\ass$ be an association scheme on $d$ classes and let
$\Ga$ be a group of algebra automorphisms of $\ass$.  Then the
matrices in $\cxa$ fixed by each element of $\Ga$ form the Bose-Mesner
algebra of a subscheme of $\ass$.
\end{theorem}

\proof
The subspace of the Bose-Mesner algebra of $\ass$ consisting of the
matrices fixed by the elements of $\Ga$ is closed under transposes,
multiplication and Schur multiplication.  Hence it is the Bose-Mesner
algebra of an association scheme, necessarily a subscheme of $\ass$.\qed

\begin{corollary}
\label{cor:hnass}
$\hnass$ is an association scheme.
\end{corollary}

\proof
Let $\Ga$ be the symmetric group on $n$ letters.  Then $\Ga$ acts as a
group of algebra automorphisms on the Bose-Mesner algebra of $\ass\pn$.
The subspace of fixed matrices is the Bose-Mesner algebra of
$\hnass$.\qed

The Bose-Mesner algebra of a scheme $\ass$ is a vector space.  Hence
we may view $\ass\pn$ as the $n$-th graded piece of the tensor
algebra.  Given this, $\hnass$ is the $n$-th graded piece of the
symmetric algebra.  Thus we have shown that the graded pieces of the
symmetric algebra are Schur-closed.  We note one useful consequence
of this.

\begin{lemma}
\label{lem:trans}
For any association scheme $\ass$ we have 
\[
    H(m,\hnass)\cong H(mn,\ass).\qed
\]
\end{lemma}

We consider one example.  Let $\clc$ denote the two-class association
scheme formed by the cycle of length four and its complement.
$\clc=H(n,2)$ and $\hn\clc=H(2n,2)$.  This observation plays a small
role in the construction of the Kerdock codes as linear codes over
$\ints_4$.  We discuss this in \sref{appls}.

\section{Generating Functions}

We start by deriving a generating function for the eigenvalues of
$\hnass$.

Let $\ass$ be an association scheme with $d$ classes, formed by the
matrices $\seq A0d$.  If $\al$ is a vector of length $d+1$ with
non-negative integer entries summing to $n$, let $A_\al$ denote the
matrix in $\hnass$ with $vw$-entry 1 if and only if $h(v,w)=\al$.
Note that $A_\al$ is the sum of all products
$$
    A_{i_1}\otimes\cdots\otimes A_{i_r}
$$
where $A_r$ occurs exactly $\al(r)$ times, for $r=0,1,\ldots,d$.  If
$\seq E0d$ are the principal idempotents of $\ass$, we similarly
define $E_\be$ to be the sum of all products
$$
    E_{i_1}\otimes\cdots\otimes E_{i_r}
$$
where $E_r$ occurs exactly $\be(r)$ times, for $r=0,1,\ldots,d$.

Let $\seq s0d$ be a set of $d+1$ independent commuting variables and
let $s$ denote the column vector with $i$-th entry equal to $s_i$ for
each $i$.  If $\al$ is a vector of length $d+1$ with non-negative
integer entries, let $s^\al$ be the monomial given by
$$
    s^\al := \prod_{i=1}^n s_i^{\al(i)}.
$$
It follows immediately that
\begin{equation}
\label{soal}
    (s_0A_0+\cdots+s_dA_d)\pn= \sum_\al s^\al A_\al,
\end{equation}
where $\al$ ranges over all $\binom{d+n}{n}$ vectors of length
$d+1$ with non-negative integer entries summing to $n$.
Similarly, if $\seq t0d$ are independent commuting variables
then
\begin{equation}
\label{tobe}
    (t_0E_0+\cdots+t_dE_d)\pn= \sum_\be t^\be E_\be.
\end{equation}

\begin{lemma}
\label{lem:atoe}
    We have
    $$
        \sum_\al A_\al s^\al =\sum_\be (Ps)^\be E_\be
    $$
    and
    $$
        \sum_\al (P^{-1} t)^\al A_\al= \sum_\be t^\be E_\be.
    $$
\end{lemma}

\proof
If
$$
    t_j:= \sum_i s_i p_i(j)
$$
then
$$
\sum_j t_jE_j= \sum_{i,j} s_ip_i(j)E_j= \sum_i s_iA_i.
$$
Hence our claims follow from \eqref{soal} and \eqref{tobe}.\qed

Suppose $\ass$ is an association scheme on $v$ vertices with matrix
of eigenvalues $P$.  An association scheme $\assb$ is {\sl formally
dual} to $\ass$ if the matrix of eigenvalues of $\assb$ is $vP^{-1}$.
Hence \lref{atoe} has the following consequence.

\begin{corollary}
\label{cor:duals}
    If $\ass$ and $\assb$ are formally dual association
    schemes then $\hnass$ and $H(n,\assb)$ are formally dual.  If $\ass$
    is formally self-dual then $\hnass$ is too.\qed
\end{corollary}

Denote the eigenvalue of $A_\al$ associated to the idempotent $E_\be$
by $p_\al(\be)$.  Then multiplying both sides of the first identity
in \lref{atoe} by $E_\ga$ yields
\begin{equation}
    \label{alsal}
    \sum_\al s^\al p_\al(\ga)= (Ps)^\ga,
\end{equation}
which provides a generating function for the eigenvalues of $\hnass$.

\section{MacWilliams Theorem}
\label{sec:macw}

Let $\ass$ be an association scheme of order $v$ with Schur
idempotents $\seq A0d$.  If $C$ is a subset of the vertices of $\ass$
with characteristic vector $x$ then the vector
$$
    \frac{1}{|C|}(x^TA_0x, x^TA_1x,\ldots,x^TA_dx)
$$
is called the {\sl inner distribution} of $C$.  The central problem of
coding theory is to find the subsets of the vertices of $\hnq$ of
maximal size, subject to given constraints on the inner distribution.
One of the main contributions of the theory of association schemes is
that it allows us to define what is essentially a Fourier transform of
the inner distribution vector.

Let $\seq E0d$ be the principal idempotents of $\ass$.  The vector
$$
    \frac{v}{|C|^2}(x^TE_0x,x^TE_1x,\ldots,x^TE_dx)
$$
is the \textsl{dual inner distribution}.  

Now suppose that $C$ is a subset of the vertices of $\hnass$,
with characteristic vector $x$.  Define the \textsl{weight enumerator}
of $C$ to be the polynomial
$$
    W_C(s) = \frac{1}{|C|}\sum_\al x^T\! A_\al x\>s^\al
$$
and the {\sl dual weight enumerator} to be
$$
    W^\perp_C(t) = \frac{v^n}{|C|^2}\sum_\be x^T\! E_\be x\>t^\be.
$$
Either of these polynomials determines the other:

\begin{theorem}
\label{thm:macw}
    Let $\ass$ be an association scheme of order $v$ and let $P$
    be its matrix of eigenvalues.  Let $C$ be a subset of the vertices of
    $\hnass$.  Then
    $$
        W^\perp_C(t) = \frac{v^n}{|C|} W_C(P^{-1} t).
    $$
\end{theorem}

\proof
We have
$$
    \sum_\al A_\al (P^{-1} t)^\al =\sum_\al E_\al t^\al,
$$
from which the result follows directly.\qed

If we restrict ourselves to translation schemes, we can say more.  An
\textsl{automorphism} of an association scheme is a permutation of its
vertices that is an automorphism of each colour class in the scheme.
(There is essentially no relation between automorphisms and
Bose-Mesner automorphisms.)  An association scheme is a 
\textsl{translation scheme} if there is an abelian group, $\Ga$ say, of
automorphisms acting transitively on its vertices.  
If $\ass$ is a
translation scheme relative to the abelian group $\Ga$, then $\hnass$
is a translation scheme relative to the abelian group $\Ga^n$.  If
$\ass$ is a translation scheme relative to $\Ga$ then, since we may
assume without loss that $\Ga$ acts faithfully and therefore regularly, 
we may identify the vertex set of $\ass$ with $\Ga$: associate some fixed vertex of
$\ass$ to the identity, and associate the image of this vertex under
an element $\ga$ of $\Ga$ with $\ga$.  A subset of the vertices of a
translation scheme is {\sl additive} if its image in $\Ga$ is a
subgroup.  A scheme may be a translation scheme relative to more than
one group; for example, $H(n,2)$ is a translation scheme relative to
the group $\ints_4^a\times\ints_2^b$ for any non-negative integers $a$
and $b$ such that $2a+b=n$.

Let $\Ga^\perp$ denote the character group of $\Ga$.  The \textsl{dual}
to a subset $C$ of $\Ga$ is the set
$$
    C^\perp := \{\chi\in\Ga^*: \chi(x)=1,\ \forall x\in C\}.
$$
Note that $C^\perp$ is a subgroup of $\Ga^\perp$ and, if $C$ is a
subgroup of $\Ga$ then $C^{\perp\perp}=C$.  As $\Ga^\perp$ and $\Ga$
are isomorphic, we may view $\Ga^\perp$ as an abelian group of
automorphisms of $\ass$, and hence identify $C^\perp$ with a subset of
the vertices of $\ass$.  With all this established \cite[Theorem~2.10.12]{bcn} 
yields:

\begin{theorem}
\label{thm:dual}
    If $\ass$ is a translation scheme and $C$ is an additive
    subset of $H(n,\ass)$, then $W_{C^\perp}(t)=W^\perp_C(t)$.\qed
\end{theorem}

This brings us to one of the main results of this paper.

\begin{corollary}
\label{cor:trmac}
    Let $\ass$ be a translation scheme and let $P$ be its
    matrix of eigenvalues.  If $C$ is an additive subset of the vertices
    of $\ass$, then
    $$
        W_{C^\perp}(t) = \frac{v^n}{|C|} W_C(P^{-1} t).\qed
    $$
\end{corollary}

We consider the simplest case, when $\ass$ is the scheme with one
class on $q$ vertices.  Then $\hnass=H(n,q)$.  A linear code of length
$n$ over $GF(q)$ is an additive subset of $H(n,q)$ and $C^\perp$ is
the usual dual code.  The matrix of eigenvalues for $\ass$ is
$$
    P=\pmat{1&q-1\\1&-1}
$$
and its inverse is $q^{-1} P$.  So \tref{macw} and \tref{dual} imply
that, if $C$ is a linear code in $H(n,q)$, then
$$
    W_{C^\perp}(x,y)=\frac1{|C|}W_C(x+(q-1)y,x-y).
$$
This is the standard form of MacWilliams theorem (see, for example,
\cite[Theorem~5.13]{MacSlo}).

\begin{corollary}
\label{cor:hompol}
    Let $\ass$ be an association scheme with $d$ classes
    and let $P$ be its matrix of eigenvalues.  Let $\widehat P$ be the
    matrix representing the action of $P$ on homogeneous polynomials
    of degree $n$ in $\seq s0d$.  Then $\widehat P$ is the matrix
    of eigenvalues of $\hnass$.\qed
\end{corollary}

Taking the Schur product of each side of the second identity in
\lref{atoe} with $A_\ga$ yields
$$
    \sum_\be t^\be q_\be(\ga)= (vP^{-1} t)^\ga,
$$
a generating function for the dual eigenvalues.

We illustrate some of these results by applying them to the Hamming
scheme.  The matrix of eigenvalues for $H(1,q)$ is
$$
    P= \pmat{1&q-1\\ 1& -1},
$$
with inverse $q^{-1} P$.  We use \lref{atoe}, obtaining
$$
\sum_{i=0}^n  s_0^{n-i} s_1^iA_{n-i,i}
	=\sum_{j=0}^n (s_0+(q-1)s_1)^{n-j}(s_0-s_1)^jE_{n-j,j};
$$
similarly
$$
q^{-n}\sum_{i=0}^n (t_0+(q-1)t_1)^{n-i}(t_0-t_1)^iA_{n-i,i}
	=\sum_{j=0}^n t_0^{n-j} t_1^j E_{n-j,j}.
$$
The generating function for the eigenvalues unfolds as:
$$
    \sum_i p_{n-i,i}(n-j,j)s_0^{n-i}s_1^i= (s_0+(q-1)s_1)^{n-j}(s_0-s_1)^j.
$$
To decode these, note that $A_{n-i,i}$ is the $i$-th minimal Schur
idempotent of $\hnq$ and $E_{n-j,j}$ is the $j$-th principal
idempotent.  Hence $p_{n-i,i}(n-j,j)$ represents $p_i(j)$.

\section{Codes}
\label{sec:appls}

We make some remarks on codes over $\ints_4$.  Our aim is to show that
a number of standard results follow easily from the theory we have
described.

By a code $C$ of length $n$ over $\ints_4$, we mean simply a subset of
$\ints_4^n$.  The relevant association scheme depends on our choice of
`distance function'.  Suppose first that, if $x$ and $y$ are words in
$C$ then $h_4(x,y)$ is the vector of length four, giving the number of
coordinates $i$ such that $x_i-y_i$ is 0, 1, 2 or 3 modulo 4.  Let
$\dcfr$ denote the three-class scheme formed by the elements of
$\ints_4$, viewed as $4\times4$ permutation matrices.  The relations
on $\ints_4^n$ determined by $h_4(x,y)$ form the scheme $\hn\dcfr$.
The matrix of eigenvalues for $\dcfr$is
$$
    P=\pmat{1&1&1&1\\ 1&i&-1&-i\\ 1&-1&1&-1\\ 1&-i&-1&i},
$$
which satisfies $P^2=4I$.  This scheme is linear and so \tref{macw}
yields
$$
W^\perp_C(s,t,u,v)
	=\frac1{|C|}W_C(s+t+u+v,s+it-u-iv,s-t+u-v,s-it-u+iv).
$$
In coding theory terms, this is MacWilliams theorem for the {\sl
complete weight enumerator} of a code in $\ints_4^n$.  Note that
$\dcfr$ is a translation scheme, so $W^\perp_C=W_{C^\perp}$ when $C$
is an additive code.  

Let $\clc$ denote the two-class association scheme belonging to the
undirected cycle on four vertices.  This is the symmetric subscheme of
$\dcfr$ and has matrix of eigenvalues
$$
    P=\pmat{1&2&1\\ 1&0&-1\\ 1&-2&1}.
$$
Now $\hn\clc$ is the scheme determined by $h_3(x,y)$, where $h_3(x,y)$
is a triple with entries counting the number of coordinates $i$ such
that $x_i-y_i$ is 0, 1 or 3, or 2 modulo 4.  If $C\sbs\ints_4^n$ then
$$
    W^\perp_C(s,t,u) =\frac1{|C|}W_C(s+2t+u,s-u,s-2t+u).
$$
This is MacWilliams theorem for the \textsl{symmetrized weight enumerator}.  
As before, $W^\perp_C=W_{C^\perp}$ when $C$ is additive.
Of course this identity is a special of case of the previous one,
which is probably best ascribed to Klemm.  (See Satz~1.5 in
\cite[Satz~1.5]{klemm}.)

The sum of the entries of $h_3(x,y)$ is the \textsl{Lee distance}
between $x$ and $y$ and
$$
    W_C(s^2,st,t^2)
$$
is a form of the \textsl{Lee weight enumerator}.

There is one surprise here.  The scheme $\clc$ is equal to $H(2,2)$
and so, by \lref{trans}, we have $H(n,\clc)\cong H(2n,2)$.  The Lee
weight enumerator for a code $C$, viewed as a subset of the vertices
of $\hn\clc$, is the usual weight enumerator for $C$, viewed as a
subset of the vertices of $H(2n,2)$, i.e., as a binary code.
However it is possible that $C$ is an additive code in $\hn\clc$,
but not in $H(2n,2)$.  In this case the binary weight enumerators for
$C$ and $C^\perp$ will satisfy
$$
    W_{C^\perp}(s^2,st,t^2)=\frac1{|C|}W_C((s+t)^2, s^2-t^2, (s-t)^2,
$$
even though $C$ and $C^\perp$ are not linear as binary codes (and
$C^\perp$ will not be the dual of $C$ in the usual sense).

For further discussion of these weight enumerators, we refer the
reader to \cite[Section~2.2]{hkcss}.

\section{Modular Invariance}
\label{sec:bbmod}

Let $\ass$ be an association scheme with $d$ classes and let $P$ be
its matrix of eigenvalues.  Following Bannai, Bannai and Jaeger
\cite{bbj}] we say that $\ass$ satisfies the \textsl{modular invariance
property} if there is a diagonal matrix $T$ and a constant $c$ such
that
$$
    (PT)^3=cI.
$$
This may seem an unlikely condition, but it does hold for the
Hamming scheme and there is an interesting connection
with the theory of spin models.

The following result generalizes Theorem~1 of Bannai and Bannai \cite{baba}.

\begin{lemma}
\label{modmod}
    If the association scheme $\ass$ satisfies the modular
    invariance property, so does $\hnass$.
\end{lemma}

\proof
Assume that $\ass$ has $d$ classes and let $P$ be its matrix of
eigenvalues.  If $M$ is a $(d+1)\times(d+1)$ matrix, let $\widehat{M}$
denote the matrix representing the induced action of $M$ on
polynomials of degree $n$.  If $T$ is diagonal then so is
$\widehat{T}$.  If $(PT)^3=cI$ then
$({\widehat{P}}\widehat{T})^3=cI$.\qed

For those whom it helps, we remark that if $\ass$ contains a matrix
$W(+)$ that determines a spin model then $W(+)\pn$ determines a spin
model and lies in $\hnass$ (not just in $\ass\pn$).

\section*{Acknowledgement}

I thank Bill Martin for a number of useful conversations on the topics in this paper.
The reader should thank him for identifying a number of mistakes in an earlier version.


\begin{thebibliography}{1}

\bibitem{baba}
Eiichi Bannai and Etsuko Bannai.
\newblock Modular invariance of the character table of the {H}amming
  association scheme {H⁢(d,⁢q)}.
\newblock {\em J. Number Theory}, 47(1):79--92, 1994.

\bibitem{bbj}
Eiichi Bannai, Etsuko Bannai, and Fran{\c{c}}ois Jaeger.
\newblock On spin models, modular invariance, and duality.
\newblock {\em J. Algebraic Combin.}, 6(3):203--228, 1997.

\bibitem{bcn}
A.~E. Brouwer, A.~M. Cohen, and A.~Neumaier.
\newblock {\em Distance-{R}egular {G}raphs}, volume~18 of {\em Ergebnisse der
  Mathematik und ihrer Grenzgebiete (3) [Results in Mathematics and Related
  Areas (3)]}.
\newblock Springer-Verlag, Berlin, 1989.

\bibitem{ACAM}
A.~{Chan} and A.~{Munemasa}.
\newblock {Hamming Graph in Nomura Algebra}.
\newblock {\em ArXiv e-prints}, October 2010.

\bibitem{delphd}
P.~Delsarte.
\newblock An algebraic approach to the association schemes of coding theory.
\newblock {\em Philips Res. Rep. Suppl.}, (10):vi+97, 1973.

\bibitem{hkcss}
A.~Roger Hammons, Jr., P.~Vijay Kumar, A.~R. Calderbank, N.~J.~A. Sloane, and
  Patrick Sol{\'e}.
\newblock The {Z4}-linearity of {K}erdock, {P}reparata, {G}oethals, and related
  codes.
\newblock {\em IEEE Trans. Inform. Theory}, 40(2):301--319, 1994.

\bibitem{klemm}
Michael Klemm.
\newblock Selbstduale {C}odes \"uber dem {R}ing der ganzen {Z}ahlen modulo {4}.
\newblock {\em Arch. Math. (Basel)}, 53(2):201--207, 1989.

\bibitem{MacSlo}
F.~J. MacWilliams and N.~J.~A. Sloane.
\newblock {\em The {T}heory of {E}rror-{C}orrecting {C}odes.}
\newblock North-Holland Publishing Co., Amsterdam, 1977.
\newblock North-Holland Mathematical Library, Vol. 16.

\bibitem{MarSt}
W.~J. Martin and D.~R. Stinson.
\newblock Association schemes for ordered orthogonal arrays and
  {(T,⁢M,⁢S)}-nets.
\newblock {\em Canad. J. Math.}, 51(2):326--346, 1999.

\end{thebibliography}
\end{document}